\documentstyle[12pt]{article}
\begin{document}
\newcommand{\equi}{\Longleftrightarrow}
\newcommand{\text}[1]{\mbox{{\rm #1}}}
\newcommand{\gd}{\delta}
\newcommand{\itms}[1]{\item[[#1]]}
\newcommand{\nin}{\in\!\!\!\!\!/}
\newcommand{\sub}{\subset}
\newcommand{\cntd}{\subseteq}
\newcommand{\go}{\omega}
\newcommand{\Pa}{P_{a^\nu,1}(U)}
\newcommand{\fx}{f(x)}
\newcommand{\fy}{f(y)}
\newcommand{\gD}{\Delta}
\newcommand{\gl}{\lambda}
\newcommand{\gL}{\Lambda}
\newcommand{\half}{\frac{1}{2}}
\newcommand{\sto}[1]{#1^{(1)}}
\newcommand{\stt}[1]{#1^{(2)}}
\newcommand{\Z}{\hbox{\sf Z\kern-0.720em\hbox{ Z}}}
\newcommand{\singcolb}[2]{\left(\begin{array}{c}#1\\#2
\end{array}\right)}
\newcommand{\ga}{\alpha}
\newcommand{\gb}{\beta}
\newcommand{\gga}{\gamma}
\newcommand{\ul}{\underline}
\newcommand{\ol}{\overline}
\newcommand{\qed}{\kern 5pt\vrule height8pt width6.5pt depth2pt}
\newcommand{\Lrraro}{\Longrightarrow}
\newcommand{\Nb}{|\!\!/}
\newcommand{\NN}{{\rm I\!N}}
\newcommand{\bsl}{\backslash}
\newcommand{\gt}{\theta}
\newcommand{\op}{\oplus}
\newcommand{\Op}{\bigoplus}
\newcommand{\CR}{{\cal R}}
\newcommand{\tr}{\bigtriangleup}
\newcommand{\grr}{\omega_1}
\newcommand{\ben}{\begin{enumerate}}
\newcommand{\een}{\end{enumerate}}
\newcommand{\ndiv}{\not\mid}
\newcommand{\bab}{\bowtie}
\newcommand{\hal}{\leftharpoonup}
\newcommand{\har}{\rightharpoonup}
\newcommand{\ot}{\otimes}
\newcommand{\OT}{\bigotimes}
\newcommand{\bwe}{\bigwedge}
\newcommand{\gep}{\varepsilon}
\newcommand{\gs}{\sigma}
\newcommand{\rbraces}[1]{\left( #1 \right)}
\newcommand{\bbox}{$\;\;\rule{2mm}{2mm}$}
\newcommand{\sbraces}[1]{\left[ #1 \right]}
\newcommand{\bbraces}[1]{\left\{ #1 \right\}}
\newcommand{\OO}{_{(1)}}
\newcommand{\TT}{_{(2)}}
\newcommand{\FF}{_{(3)}}
\newcommand{\minus}{^{-1}}
\newcommand{\CV}{\cal V}
\newcommand{\CVs}{\cal{V}_s}
\newcommand{\un}{U_q(sl_n)'}
\newcommand{\on}{O_q(SL_n)'}
\newcommand{\slq}{U_q(sl_2)}
\newcommand{\olq}{O_q(SL_2)}
\newcommand{\UU}{U_{(N,\nu,\go)}}
\newcommand{\HH}{H_{n,q,N,\nu}}
\newcommand{\ct}{\centerline}
\newcommand{\bs}{\bigskip}
\newcommand{\qua}{\rm quasitriangular}
\newcommand{\ms}{\medskip}
\newcommand{\noin}{\noindent}
\newcommand{\mat}[1]{$\;{#1}\;$}
\newcommand{\raro}{\rightarrow}
\newcommand{\map}[3]{{#1}\::\:{#2}\raro{#3}}
\newcommand{\C}{{\bf C}}
\newcommand{\alg}{{\rm Alg}}
\def\newtheorems{\newtheorem{theorem}{Theorem}[subsection]
                 \newtheorem{cor}[theorem]{Corollary}
                 \newtheorem{prop}[theorem]{Proposition}
                 \newtheorem{lemma}[theorem]{Lemma}
                 \newtheorem{defn}[theorem]{Definition}
                 \newtheorem{Theorem}{Theorem}[section]
                 \newtheorem{Corollary}[Theorem]{Corollary}
                 \newtheorem{Proposition}[Theorem]{Proposition}
                 \newtheorem{Lemma}[Theorem]{Lemma}
                 \newtheorem{Definition}[Theorem]{Definition}
                 \newtheorem{Example}[Theorem]{Example}
                 \newtheorem{Remark}[Theorem]{Remark}
                 \newtheorem{claim}[theorem]{Claim}
                 \newtheorem{sublemma}[theorem]{Sublemma}
                 \newtheorem{example}[theorem]{Example}
                 \newtheorem{remark}[theorem]{Remark}
                 \newtheorem{question}[theorem]{Question}
                 \newtheorem{conjecture}{Conjecture}[subsection]
                 \newtheorem{Conjecture}[Theorem]{Conjecture}
                 \newtheorem{Question}[Theorem]{Question}}
\newtheorems
\newcommand{\proof}{\par\noindent{\bf Proof:}\quad}
\newcommand{\dmatr}[2]{\left(\begin{array}{c}{#1}\\
                            {#2}\end{array}\right)}
\newcommand{\doubcolb}[4]{\left(\begin{array}{cc}#1&#2\\
#3&#4\end{array}\right)}
\newcommand{\qmatrl}[4]{\left(\begin{array}{ll}{#1}&{#2}\\
                            {#3}&{#4}\end{array}\right)}
\newcommand{\qmatrc}[4]{\left(\begin{array}{cc}{#1}&{#2}\\
                            {#3}&{#4}\end{array}\right)}
\newcommand{\qmatrr}[4]{\left(\begin{array}{rr}{#1}&{#2}\\
                            {#3}&{#4}\end{array}\right)}
\newcommand{\smatr}[2]{\left(\begin{array}{c}{#1}\\
                            \vdots\\{#2}\end{array}\right)}

\newcommand{\ddet}[2]{\left[\begin{array}{c}{#1}\\
                           {#2}\end{array}\right]}
\newcommand{\qdetl}[4]{\left[\begin{array}{ll}{#1}&{#2}\\
                           {#3}&{#4}\end{array}\right]}
\newcommand{\qdetc}[4]{\left[\begin{array}{cc}{#1}&{#2}\\
                           {#3}&{#4}\end{array}\right]}
\newcommand{\qdetr}[4]{\left[\begin{array}{rr}{#1}&{#2}\\
                           {#3}&{#4}\end{array}\right]}

\newcommand{\qbracl}[4]{\left\{\begin{array}{ll}{#1}&{#2}\\
                           {#3}&{#4}\end{array}\right.}
\newcommand{\qbracr}[4]{\left.\begin{array}{ll}{#1}&{#2}\\
                           {#3}&{#4}\end{array}\right\}}

\title{On the Exponent of Finite-Dimensional Hopf Algebras}
\author{Pavel Etingof
\\Department of Mathematics\\
Harvard University\\Cambridge, MA 02138
\and Shlomo Gelaki\\Department of Mathematics\\
University of Southern California\\Los Angeles, CA 90089}
\date{December 28, 1998}
\maketitle

\section{Introduction}
One of the classical notions of group theory is the notion
of the exponent of a group. The exponent of a group is the least common
multiple of orders of its elements.

In this paper we generalize the notion of exponent to Hopf algebras.
We define the exponent of a Hopf algebra $H$ to be the smallest $n$ such
that $m_n\circ (I\otimes S^{-2}\otimes\cdots\otimes S^{-2n+2})\circ
\Delta_n=
\varepsilon\cdot 1,$ where $m_n,$ $\Delta_n,$ $S,$ $1,$ and $\varepsilon$
are the
iterated product and coproduct, the antipode, the unit, and the counit.
If $H$ is involutive (for example, $H$ is semisimple and cosemisimple),
the last formula reduces to $m_n\circ \Delta_n=\varepsilon\cdot 1$.

We give four other equivalent definitions of the exponent (valid for
finite-dimensional Hopf algebras). In particular, we show
that the exponent of $H$ equals the order of the Drinfeld element $u$
of the quantum double $D(H),$ and the order of ${\cal R}_{21}{\cal R},$
where ${\cal R}$ is the universal $R$-matrix of $D(H).$

We show that the exponent is invariant under twisting. We prove that for
semisimple and cosemisimple Hopf algebras $H,$ the exponent is finite
and divides $dim(H)^3.$ For triangular semisimple Hopf algebras in
characteristic zero, we show that the exponent divides $dim(H).$ These
theorems are motivated by the work of Kashina [Ka1,Ka2], who
conjectured that if $H$ is semisimple and cosemisimple then
(using our language) the exponent of $H$ always divides $dim(H),$ and
showed that the order of the squared antipode of any
finite-dimensional semisimple and cosemisimple Hopf algebra in the
Yetter-Drinfeld category of $H$ divides the exponent of $H.$

At the end we formulate some open questions, in particular
suggest a formulation for a possible Hopf algebraic analogue of
Sylow's theorem.

\noindent
{\bf Acknowledgment} The authors thank Yevgenia Kashina for suggesting the
problem, and Susan Montgomery for useful discussions.

\section{Definition and Elementary Properties of Exponent}
Let $H$ be a Hopf algebra over any field $k,$ with
multiplication map $m,$ comultiplication map $\Delta$ and antipode $S.$
Let $m_1=I$ and $\Delta_1=I$ be the identity map $H\to H$, and for any
integer $n\ge 2$ let $m_n:H^{\ot n}\raro H$ and
$\Delta_n:H\raro H^{\ot n}$ be defined by $m_n=m\circ(m_{n-1}\ot
I),$ and $\Delta_n=(\Delta_{n-1}\ot I)\circ\Delta$.
We start by making the following definition.
\begin{Definition}\label{exp}
The exponent of a Hopf algebra $H,$ denoted by $exp(H),$ is the
smallest positive integer $n$ satisfying $m_n\circ (I\ot S^{-2}\ot \cdots
\ot S^{-2n+2})\circ \Delta_n=\varepsilon\cdot 1.$ If such $n$ does not
exist, we say that $exp(H)=\infty.$
\end{Definition}

Let us list some of the elementary properties of $exp(H).$
\begin{Proposition}\label{elem}
Let $H$ be a Hopf algebra over $k.$ Then:
\ben
\item
The order of any grouplike element of $H$ divides $exp(H)$ (here we agree
that any positive integer $n$ divides $\infty$).
\item
For any group $G,$ $exp(k[G])$ equals the exponent
of G (see e.g. [Ro, p. 12]), i.e. the least common multiple of
the orders of the elements of $G.$
\item
If $exp(H)=n<\infty$ then
$m_{r}\circ(I\ot S^{-2}\ot \cdots
\ot S^{-2r+2})\circ\Delta_{r}=\varepsilon\cdot 1$ if and
only if $r$ is divisible by $n.$
\item
If $H$ is finite-dimensional, $exp(H^*)=exp(H).$
\item
$exp(H_1\ot H_2)$ equals the least common multiple of $exp(H_1)$
and $exp(H_2).$
\item
If $exp(H)=2$ then $H$ is commutative and cocommutative (this
generalizes the fact that a group $G$ with $g^2=1$
for all $g\in G$ is abelian).
\item
The exponents of Hopf subalgebras and quotients of $H$ divide $exp(H).$
\item
If $K\supseteq k$ is a field then $exp(H\ot _k K)=exp(H).$
\een
\end{Proposition}
\proof 1. Suppose $exp(H)<\infty,$ and set $n=exp(H).$
Since $S^2(g)=g$ we have that $g^n=m_n\circ \Delta_n(g)=m_n\circ (I\ot
S^{-2}\ot \cdots \ot S^{-2n+2})\circ \Delta_n(g)=\varepsilon(g)1=1.$
Therefore the order of $g$ divides $n.$

\noindent
6. Since $m\circ(I\ot S^{-2})\circ\Delta=\varepsilon\cdot 1$ is
equivalent to
$S^3=I,$ we have that $I:H\raro H$ is an antiautomorphism of algebras and
coalgebras, and the result follows.

\noindent
The proofs of the other parts are obvious.
\qed
\begin{Remark}\label{motiv}
{\rm Part 2 of Proposition \ref{elem} motivated Definition \ref{exp}.}
\end{Remark}
\begin{Example}\label{inford}
Let $H$ be a finite-dimensional Hopf algebra over an algebraically closed
field $k$ of characteristic zero. Suppose that $H$ contains
a non-trivial $1:g$ skew-primitive element $x$ (i.e.
$\Delta(x)=x\ot 1 +g\ot x,$ where $g$ is a grouplike element, and
$x\notin k[g]$). It is clear that in this case we may assume that $xg=qgx$
for some root of unity $q$ of order dividing $|g|.$ Also, $S^2(x)=qx,$
$\varepsilon(x)=0,$ $\{x,gx,\dots,g^{|g|-1}x\}$ is linearly independent,
and hence $m_n\circ (I\ot S^{-2}\ot
\cdots \ot S^{-2n+2})\circ \Delta_n(x)=\sum_{i=0}^{n-1}q^{-i}g^ix\ne 0.$
Hence, $exp(H)=\infty.$ In particular, the exponent of any pointed Hopf
algebra $H$ over $k$ (which is not a group algebra)
is $\infty,$ since by [TW], $H$ contains a
non-trivial skew-primitive element.
\end{Example}

In the sequel, we will assume for simplicity that $H$ is
finite-dimensional. Let us formulate four equivalent
definitions of $exp(H).$ Recall that the Drinfeld double $D(H)=H^{*cop}\ot
H$ of $H$ is a quasitriangular Hopf algebra with universal $R-$matrix
${\cal R}=\sum_ih_i\ot h_i^*,$ where $\{h_i\},$ $\{h_i^*\}$ are dual
bases for
$H$ and $H^*$ respectively. Let $u=m(S\ot I)\tau({\cal R})=\sum_iS(h_i^*)
h_i,$ where $S$ is the antipode of $D(H)$ and $\tau$ is the usual flip
map, be the Drinfeld element of $D(H).$ By [D],
$$S^2(x)=uxu^{-1},\;x\in D(H)\;\;and\;\;\Delta(u)=(u\ot u)({\cal
R}_{21}{\cal R})^{-1}.$$
\begin{Theorem}\label{Thm1}
Let $H$ be a finite-dimensional Hopf algebra over $k.$ Then
\ben
\item
$exp(H)$ equals the smallest positive integer $n$ such that
$${\cal R}(I \ot S^{2})({\cal R})\cdots (I\ot S^{2n-2})({\cal R})=1.$$
\item
$exp(H)$ equals the order of $u.$
\item
$exp(H)$ equals the order of ${\cal R}_{21}{\cal R}.$
\item
$exp(H)$ equals the order of any non-zero element $v\in D(H)$
satisfying $$\Delta(v)=(v\ot v)({\cal R}_{21}{\cal R})^{-1}.$$
\een
\end{Theorem}
\proof
First note that since $(\Delta\ot I)({\cal R})={\cal R}_{13}{\cal
R}_{23},$ it follows that
$$(\Delta_n\ot I)({\cal R})={\cal R}_{1,n+1}\cdots{\cal R}_{n,n+1}$$
for all $n.$

Second, recall that the map $H^{*}\ot H\raro D(H),$ $p\ot h\mapsto ph$ is
a linear isomorphism [D].

Now we will show the equivalence of Definition 2.1 and the four definitions
in the theorem.

\noindent
${\rm (Definition\;\ref{exp}\Leftrightarrow 1\Leftrightarrow 2)}$ Write
${\cal R}=\sum_ja_{j}\ot b_{j}.$ Using the above we obtain the following
equivalences: \[
\begin{array}{ll}
m_n\circ (I\ot S^{-2}\ot \cdots \ot S^{-2n+2})\circ
\Delta_n=\varepsilon\cdot 1 & \equi \\ \cr
(m_n\circ (I\ot S^{-2}\ot \cdots \ot S^{-2n+2})\circ
\Delta_n\ot I)({\cal R})=1 & \equi \\ \cr
(m_n\circ (I\ot S^{-2}\ot \cdots \ot S^{-2n+2})\ot I)
({\cal R}_{1,n+1}\cdots{\cal R}_{n,n+1})=1 & \equi \\ \cr
\sum_{i_1,\dots,i_n} a_{i_1}S^{-2}(a_{i_2})\cdots
S^{-2n+2}(a_{i_n})\ot b_{i_1}\cdots b_{i_n}=1 & \equi \\ \cr
{\cal R}(I \ot S^{2})({\cal R})\cdots (I\ot S^{2n-2})({\cal R})=1 & \equi
\\ \cr
u^n=1
\end{array}
\]
(in the last step we applied $m\circ(I\ot S)\tau$ to both sides of the
equation, and used the fact that $uS^{-2}(x)=xu,$ for all $x\in D(H)$).

\noindent
${\rm (2\Leftrightarrow 3)}$ Clearly if $u^n=1$ then $({\cal R}_{21}{\cal
R})^n=1.$ In the other direction, first note that $({\cal
R}_{21}{\cal R})^n=1$ implies that $u^n\in G(D(H))$
(where $G(A)$ is the group of grouplike elements of a Hopf algebra
$A$). Therefore by [R], $u^n=ab$ where $a\in G(H^*)$ and $b\in G(H).$
Regarding $u$ as an element of $H^*\ot H,$ we have that
$m(I\ot \varepsilon)(u)=m(\varepsilon\ot I)(u)=1.$ Hence it follows that
$1=m(I\ot \varepsilon)(u^n)=a$ and $1=m(\varepsilon\ot I)(u^n)=b,$ so
$u^n=1.$

\noindent
${\rm (2\Leftrightarrow 4)}$ First note that $v=ug,$ where $g\in G(D(H)).$
Since $g$ commutes with $u$ we have that $v^n=u^ng^n.$ Therefore if
$u^n=1$ then $v^n=1$ by parts 1 and 3 of Proposition \ref{elem}, and if
$v^n=1$ then $u^n\in G(D(H)),$ so $u^n=1$ as explained above. \qed
\begin{Corollary}\label{dual}
Let $H$ be a finite-dimensional Hopf algebra over $k.$ Then
$$exp(H^{cop})=exp(H^{op})=exp(H).$$
\end{Corollary}
\proof Since $(D(H^{*cop}),{\tilde{\cal
R}})\cong (D(H)^{op},{\cal R}_{21})$ as quasitriangular Hopf algebras,
it follows from part 1 of Theorem \ref{Thm1} that $exp(H^{*cop})=exp(H).$
Hence the result follows from part 4 of Proposition \ref{elem}. \qed
\section{Invariance of Exponent Under Twisting}
In this section we show that $exp(H)$ is invariant under twisting.

First recall Drinfeld's notion of a twist for Hopf algebras.
\begin{Definition}\label{twist}
Let $H$ be a Hopf algebra over $k.$ A twist for $H$
is an invertible element $J\in H\ot H$ which satisfies
$$
(\Delta\ot I)(J)J_{12}=(I\ot \Delta)(J)J_{23}\;\; and \;\; (\varepsilon\ot
I)(J)=(I\ot \varepsilon)(J)=1.
$$
\end{Definition}
Given a twist $J$ for $H$, we can construct a new Hopf algebra $H^J$,
which is the same as $H$ as an algebra, with coproduct $\Delta^J$ given by
$$\Delta^J(x)=J^{-1}\Delta(x)J,\; x\in H.$$
If $(H,R)$ is quasitriangular then so is $H^J$ with the $R$-matrix
$$R^J=J_{21}^{-1}RJ.$$
In particular, since $H$ is a Hopf subalgebra of $D(H),$ we can twist
$D(H)$ using the twist $J\in D(H)\ot D(H)$ and obtain $(D(H)^J,{\cal R}^J).$
\begin{Proposition}\label{tdouble}
Let $H$ be a finite-dimensional Hopf algebra over $k,$  and let $J$ be a
twist for $H.$ Then
$(D(H^J),{\cal R})\cong (D(H)^J,{\cal R}^J)$ as quasitriangular Hopf
algebras.
\end{Proposition}
\proof Let $H_+$ and $H_-$ be the Hopf subalgebras of $D(H)^J$ generated
by the left and right components of ${\cal R}^J$ respectively. Clearly,
$H_+\subseteq H^J.$ In order to prove the theorem it is sufficient to
prove that the multiplication map $H_+\ot H_-\raro D(H)^J$ is a linear
isomorphism, since then $H_+=H^J$ (by dimension counting) and the result
will follow.

Clearly, $dim(H_+)\le dim(H)$ and $dim(H_-)=dim(H_+),$ so we need to show
that $H_+H_-=D(H).$ Since $J{\cal R}_{21}^J{\cal
R}^JJ^{-1}={\cal R}_{21}{\cal R}$ we have that $HH_-H=D(H)$ (looking at
the first component). Let
$A=H_+H_-=H_-H_+,$ $dim(H)=d,$ $dim(H_+)=d_+,$ $\{v_1,\dots,v_{d/d_+}\}$
with $v_1=1$ be a
basis of $H$ as a right $H_+-$module, and $\{w_1,\dots,w_{d/d_+}\}$ with
$w_1=1$ be a
basis of $H$ as a left $H_+-$module (such bases exist by the freeness
theorem for Hopf algebras [NZ]). Then we get by dimension counting that
$D(H)=\bigoplus_{i,j=1}^{d/d_+}v_iAw_j.$ Thus, $HH_-\cap H_-H=A,$
hence $H\subseteq A$ which implies that $HAH=A,$ and the result follows.
\qed
\begin{Theorem}\label{tinv}
Let $H$ be a finite-dimensional Hopf algebra over $k,$ and let $J$ be a
twist for $H.$ Then $exp(H)=exp(H^J).$
\end{Theorem}
\proof
By part 3 of Theorem \ref{Thm1}, and Proposition \ref{tdouble}, it is
sufficient to show that the order of ${\cal R}_{21}^J{\cal R}^J$ equals
to the order of ${\cal R}_{21}{\cal R}.$ But this is clear since they are
conjugate. \qed
\begin{Corollary}\label{ddouble}
Let $H$ be a finite-dimensional Hopf algebra over $k.$ Then
$exp(D(H))=exp(H).$
\end{Corollary}
\proof
By [RS], there exists $J\in D(H)\ot D(H)$ such that
$D(D(H))\cong (D(H)\ot D(H))^J$
as Hopf algebras. Then using Theorem \ref{tinv} we get that
$exp(D(H))$ equals the order of $u$ in
$(D(H)\ot D(H))^J$ which equals the order of $u$ in $D(H)\ot D(H),$
and hence equals $exp(H)$ (since $\displaystyle{u_{_{D(H)\ot
D(H)}}=u_{_{D(H)}}\ot u_{_{D(H)}}}$). \qed
\section{The Exponent of a Semisimple and Cosemisimple Hopf Algebra}
In this section, we will show that if $H$ is semisimple and cosemisimple
then $exp(H)$ is finite, and give an estimate for it in terms of $dim(H).$

Let $H$ be a semisimple and cosemisimple Hopf algebra over $k$
(note that by [LR] the cosemisimplicity assumption is redundant if the
characteristic of $k$ is $0$). Recall that for semisimple and
cosemisimple $H,$ $D(H)$ is also semisimple and cosemisimple [R]. Also, by
[LR, EG2], $S^2=I$ and hence $u$ is central in $D(H).$ This implies that
$exp(H)$ equals the smallest positive integer $n$ satisfying
$m_n\circ \Delta_n=\varepsilon\cdot 1,$
and also to the order of $\mathcal R$ (by part 1 of Theorem
\ref{Thm1}).
\begin{Remark}\label{kashina}
{\rm In [Ka1,Ka2] Kashina studied the smallest positive integer $n$
satisfying
$m_n\circ \Delta_n=\varepsilon\cdot 1,$ for any finite-dimensional
Hopf algebra $H.$ In particular she observed the analogous properties
listed in Proposition \ref{elem}, and proved an analogue to Corollary
\ref{ddouble} under the assumption that this smallest $n$ is the same for
$H$ and $H^{cop}.$}
\end{Remark}
\begin{Theorem}\label{t}
Let $(H,R)$ be a semisimple triangular Hopf algebra over a
field $k$ of characteristic $0.$ Then $exp(H)$ divides $dim(H).$
\end{Theorem}
\proof By part 8 of Proposition \ref{elem}, we may assume that $k$ is
algebraically closed. Now,
it is straightforward to check that the theorem holds for
$(k[G],1\ot 1)$ where $G$ is a finite group. But by [EG1, Theorem 2.1],
there
exist a finite group $G$ and a twist $J\in k[G]\ot k[G]$ such that
$H\cong k[G]^J$ as Hopf algebras. Hence the result follows from
Theorem  \ref{tinv}. \qed
\begin{Theorem}\label{main}
Let $H$ be a semisimple and cosemisimple Hopf algebra over $k.$ Then
$exp(H)$ divides $dim(H)^3.$
\end{Theorem}
Theorem \ref{main} will be proved later.
\begin{Corollary}\label{yd}
Let $H$ be a semisimple and cosemisimple Hopf algebra over $k,$ and let
$B$ be a finite-dimensional semisimple Hopf algebra in the category of
Yetter-Drinfeld modules over $H.$ Then the order of the antipode of $B$
is finite and divides $2dim(H)^3,$ and if $H$ is semisimple triangular and
the characteristic of $k$ is $0,$ then the order of the antipode of $B$
is finite and divides $2dim(H).$
\end{Corollary}
\proof
Follows from Theorems \ref{t} and \ref{main}, and [Ka1, Theorem 6]. \qed
\begin{Remark}\label{vafa}
{\rm Theorem \ref{main} is motivated by Vafa's theorem [V].
Vafa's theorem (see [Ki] for the mathematical exposition) states
that the twists in a semisimple modular category act on the irreducible
objects by multiplication by roots of unity. Thus,
the fact that $u\in D(H)$ has a finite order follows from the fact that the
category of representations of $D(H)$ is modular, with
system of twists given by the action of the central element $u$ (see e.g.
[EG1]).}
\end{Remark}

Kashina conjectured the following:

\begin{Conjecture}\label{dimh}
Let $H$ be a semisimple and cosemisimple Hopf algebra over $k.$ Then
$exp(H)$ divides $dim(H).$
\end{Conjecture}

This conjecture was checked by Kashina in a number of special cases
[Ka1,Ka2]. Our results presented above give additional supportive evidence
for this conjecture.

Now we will prove Theorem \ref{main}. In order to do this, we need a
lemma.
\begin{Lemma}\label{thm1}
Let $H$ be a Hopf algebra of finite dimension $d$ over $k,$ ${\cal R}\in
H\ot H^{*cop}\subset D(H)\ot D(H)$ be the universal R-matrix, and $u\in
D(H)$ be the Drinfeld element. Then:
\ben
\item
For any finite-dimensional $H-$module $V_+$ and finite-dimensional
$H^*-$module $V_-,$ one has
$(det({\cal R}_{|V_+\ot V_-}))^d=1.$
\item
For any finite-dimensional $D(H)-$module $V,$ one has
$(det(u_{_{|V}}))^{d^2}=1.$
\een
\end{Lemma}
\proof 1. Recall that $(\Delta\ot I)({\cal R})={\cal R}_{13}{\cal
R}_{23}.$ Apply this identity to $V_+\ot
H\ot V_-,$ where $H$ is the regular representation of $H.$ Since
$V_+\ot H=(dimV_+)H,$ this yields, after taking determinants:
$$(det({\cal R}_{|H\ot V_-}))^{dimV_+}=(det({\cal R}_{|V_+\ot
V_-}))^d(det({\cal R}_{|H\ot
V_-}))^{dimV_+}.$$ The result follows after cancellation.

\noindent
2. We use Drinfeld's formula, $\Delta(u)=(u\ot u)({\cal
R}_{21}{\cal R})^{-1}.$ Using part 1, and the fact that $D(H)=H^*\ot H$
as $H^*-$module and $H-$module, we compute:
$$det(\Delta(u)_{|V\ot
D(H)})=(det(u_{_{|V}}))^{d^2}(det(u_{_{|D(H)}}))^{dimV}.$$
Since $V\ot D(H)=(dimV)D(H),$ the result follows. \qed

\noindent
{\bf Proof of Theorem \ref{main}:}
By part 8 of Proposition \ref{elem}, we may assume that $k$ is
algebraically closed. Since $u$ is central, we have for any irreducible
$D(H)-$module $V$ that
$det(u_{_{|V}})=\lambda(u,V)^{dimV},$ where $\lambda(u,V)$ is the
eigenvalue of $u$ on $V.$ So by
Lemma \ref{thm1}, $\lambda(u,V)^{dimV\cdot d^2}=1.$ But by
[EG1, Theorem 1.5], and in positive characteristic by [EG2, Theorem 3.7],
$dimV$ divides $d,$ so $\lambda(u,V)^{d^3}=1.$
Thus, $u^{d^3}=1$ and we are done by part 2 of Theorem \ref{Thm1}. \qed

In the non-semisimple case, as we know, Theorem \ref{main} fails, and the
order of $u$ may be infinite. The analogue of Theorem \ref{main} in this
case is the following theorem.

Let $A$ be a finite-dimensional algebra.
For any two irreducible $A-$modules $V_1$ and $V_2,$ write $V_1\sim V_2$ if
they occur as constituents in the same indecomposable $A-$module. Extend
$\sim$ to an equivalence relation. For an irreducible module $W,$ let
$[W]$ be the equivalence class of $W.$ For an indecomposable module $V$
let $[V]$ be the equivalence class of any constituent $W$ of $V.$ Let
$N_{_V}$ be the greatest common divisor of dimensions of elements of
$[V].$
\begin{Theorem}\label{last}
Let $H$ be a Hopf algebra of dimension $d$ over $k.$ Then:
\ben
\item
For any indecomposable $D(H)-$module $V,$ the unique eigenvalue of the
central element $z=uS(u)$ on $V$ is a root of unity of order dividing
$d^2N_{_V}.$
\item
For any indecomposable $D(H)-$module $V,$ every eigenvalue of $u$ on $V$
is a root of unity of order dividing $2d^2N_{_V}$ (so the
eigenvalues of $u$ on any $D(H)$ module are roots of unity).
\een
\end{Theorem}
\proof By part 8 of Proposition \ref{elem}, we may assume that $k$ is
algebraically closed.

\noindent
1. Recall that $z=u^2g,$ where $g$ is a grouplike element of
$D(H).$ By [NZ], the order of $g$ divides $dim(H)=d.$
Thus we have $(det(u_{_{|V}}))^{2d}=(det(z_{_{|V}}))^{d}.$ Since $V$
is indecomposable and $z$ is central, $z$ has a unique eigenvalue
$\lambda(z,V)$ on $V.$ For any $W\in [V],$ $\lambda(z,V)=\lambda(z,W),$
so we get,
$(det(u_{_{|W}}))^{2d}=\lambda(z,V)^{d\cdot dimW},$ which implies by
part 2 of Lemma \ref{thm1}, that $\lambda(z,V)^{d^2\cdot dimW}=1.$

\noindent
2. Since, any eigenvalue
$\mu$ of $u_{_{|V}}$ has the property $\mu ^2=\lambda(z,V)\nu,$ where
$\nu$ is an eigenvalue of $g^{-1},$ we have by part 1 that $\mu
^{2d^2dimW}=1,$
and the result follows. \qed
\begin{Corollary}\label{nss}
If $exp(H)=\infty$ then $u$ is not semisimple.
\end{Corollary}
\begin{Corollary}\label{fin}
Let $H$ be a finite-dimensional Hopf algebra over a field $k$ of
positive characteristic $p.$ Then $exp(H)<\infty.$
\end{Corollary}
\proof Let $u$ be the Drinfeld element of $D(H).$ By part 2
of Theorem \ref{last}, the eigenvalues of $u$ are roots of unity. Hence
there exists a positive integer $a$ such that $u^a=1+n$ where $n\in
D(H)$ is a nilpotent element. But then $u^{ap^b}=1$ for a sufficiently
large positive integer $b,$ and the result follows from part 2 of Theorem
\ref{Thm1}. \qed
\section{Concluding Remarks}
In conclusion we would like to formulate some questions.

\begin{Question}\label{1} Suppose that $H$ is a semisimple
and cosemisimple Hopf algebra of dimension $d$ over $k.$
If a prime $p$ divides $d,$ must it divide $exp(H)$?
\end{Question}
We do not know the answer to this question even in characteristic zero,
even for $p=2.$ However, if $H$ is a group algebra then the answer is
positive, since the statement is equivalent to (a special case of)
Sylow's first theorem: a finite group whose order is divisible by $p$ has
an element of order $p.$ So positive answer to Question \ref{1} would be a
"quantum Sylow theorem".
\begin{Question}\label{2} Let $H$ be a semisimple
and cosemisimple Hopf algebra over $k$ whose exponent is a power of a
prime $p.$
Must the dimension of $H$ be a power of the same prime?
\end{Question}
This is a special case of Question \ref{1}, but we still do not know the
answer, except for the case $exp(H)=2,$ when the answer is trivially
positive. For group algebras, the statement is equivalent to the well-known
group-theoretical result that a finite group where orders of all elements
are powers of $p$ is a p-group (a special case of Sylow's theorem).
\begin{Question}\label{3} Let $H$ be a finite-dimensional Hopf algebra
over $k$ such
that the element $u\in D(H)$ is semisimple in the regular representation.
Does it follow that $H$ and $H^*$ are semisimple
\ben
\item
In characteristic zero?
\item
In positive characteristic $p$?
\een
\end{Question}
By Theorem \ref{last}, Part 1 of Question \ref{3} is equivalent to the
question whether
for a finite-dimensional Hopf algebra $H$ in characteristic $0,$
$exp(H)<\infty$ implies that $H$ is semisimple.

A positive answer to part 2 of Question \ref{3} implies a positive
answer to Question \ref{1} for involutive Hopf algebras defined over $\Z$
and free as $\Z$-modules (which includes group algebras, i.e.
this would generalize Sylow's theorem). Indeed, if $H$ is such a Hopf
algebra then for any prime $p$ dividing the dimension of $H,$
either $H/pH$ or $(H/pH)^*$ is not semisimple (as $tr(S^2)=0$),
and hence $D(H)$ is not semisimple. If the answer to part 2 of Question
\ref{3} is positive, then this would imply that $u$ is not semisimple
over $F_p,$ i.e. the order of $u$ is divisible by $p,$ as desired.

For group algebras, the answer to part 2 of Question \ref{3} is positive: in
this case semisimplicity of $u$ is equivalent to semisimplicity of
$R=\sum g\ot \delta_g$, which implies that all group elements $g$ are
semisimple. This would imply that their orders are not divisible by $p,$
which by Sylow's theorem implies that $p$ does not divide the order of the
group.

\end{document}